\documentclass[b5paper,reqno]{amsart}
\usepackage{amsfonts,lingmacros,tree-dvips, amsmath, amssymb, graphicx, mathrsfs}
\usepackage[hidelinks]{hyperref}

\usepackage[all]{xy}

\newdir{ >}{!/6pt/@{ }*:(1,0.2)@_{>}*:(1,-0.2)@^{>}}

\theoremstyle{plain}

\newtheorem{theorem}{Theorem}[section]
\newtheorem{lemma}[theorem]{Lemma}

\theoremstyle{definition}
\newtheorem{definition}[theorem]{Definition}

\newtheorem{remark}[theorem]{Remark}
\newtheorem{counter example}[theorem]{Counter Example}

\numberwithin{equation}{section}

\DeclareMathAlphabet{\mathscr}{OT1}{pzc}{m}{it} 
\begin{document}
\Large{
		\title{ ON NON-BAIRE RARE SETS IN CATEGORY BASES}
		
		\author[S. Basu]{Sanjib Basu}
		\address{\large{Department of Mathematics,Bethune College,181 Bidhan Sarani}}
		\email{\large{sanjibbasu08@gmail.com}}
		
		\author[A.C.Pramanik]{Abhit Chandra Pramanik}
	\address{\large{Department of Pure Mathematics, University of Calcutta, 35, Ballygunge Circular Road, Kolkata 700019, West Bengal, India}}
	\email{\large{abhit.pramanik@gmail.com}}

		\thanks{The second author thanks the CSIR, New Delhi – 110001, India, for financial support}
	\begin{abstract}
  In this paper, we deal with non-Baire rare sets in category bases which forms $\aleph_0$-independent family, where a rare set is a common generalization of both Luzin and Sierpinski set. 
	\end{abstract}
\thanks{}
	\maketitle

\section{INTRODUCTION}
In $1914$, Luzin [$4$] constructed using continuum hypothesis ($2^{\aleph_0}=\aleph_1$) an uncountable set of reals having countable intersection with every set of first category. A similar construction was given by Mahlo [$4$] a year before in $1913$, but in literature this set is commonly known as the Luzin set probably because Luzin investigated these types of sets more thoroughly and proved a number of its important properties.\\
The dual of a Luzin set is the Sierpinski set constructed by Sierpinski [$4$] in $1924$ using the same continuum hypothesis. It is an uncountale set of reals having countable intersection with every set of Lebesgue measure zero. From the dual nature of $\sigma$-ideal of first category sets and the $\sigma$-ideal of Lebesgue null sets and their complementarity in the real line, it follows that the nature of the above two sets are dual of each other. For example, a Luzin set is a set of measure zero (in fact, they are universally null sets) [$4$], whereas a Sierpinski set is a set of first category; no uncountable subset of a Luzin set can have the Baire property, whereas no uncountable subset of a Sierpinski set is Lebesgue measurable etc [$4$] and they follow from some instrisic properties of the real line that are closely connected with measure and category. This paper involves the notion of `rare set' in a category base which is a common generalization of both a Luzin set and its dual the Sierpinski set. 

\section{PRELIMINARIES AND RESULTS}

The concept of a category base is a generalization of both measure and topology. Its main objective is to present both measure and Baire category (topology) and also some other aspects of point set classification within a common framework. It was introduced by J. C. Morgan II in the seventies of the last century and since then has been developed through a series of papers [$1$], [$8$], [$9$], [$11$] etc.
\begin{definition}[\cite{QM1979}]
A pair $(X,\mathcal{C}$) where $X$ is a non-empty set and $\mathcal{C}$ is a family of subsets of $X$ is called a category base if the non-empty members of $\mathcal{C}$ called regions satisfy the following set of axioms:
\begin{enumerate}
\item Every point of $X$ belongs to some region; i.e., $X = \cup$ $\mathcal{C}$.
\item Let $A$ be a region and $\mathcal{D}$ be any non-empty family of disjont regions having cardinality less than the cardinality of $\mathcal{C}$.\\
i) If $A \cap ( \cup \mathcal{D}$) contains a region, then there exists a region $D\in\mathcal{D}$ such that $A\cap D$ contains a region. \\
ii)  If $A\cap(\cup \mathcal{D})$ contains no region, then there exists a region $B\subseteq A$ which is disjoint from every region in $\mathcal{D}$.
\end{enumerate}
\end{definition}

\begin{definition}[\cite{QM1979}]
In a category base ($X,\mathcal{C}$), a set is called `singular' if every region contains a subregion which is disjoint from the set itself. A set which can be expressed as countable union of singular sets is called `meager'. Otherwise, it is called `abundant'. In particular, it is abundant everywhere (resp. everywhere in a region) if the set is abundant in every region (resp. in every subregion of it). A set is `Baire' (or, a Baire set) if every region contains a subregion in which the set or its complement is meager. 
\end{definition}
\begin{definition}[\cite{QM1979}]
A category base ($X,\mathcal{C}$) is called `point-meager' if every singleton set in it is meager and a `Baire base' if every region in it is abundant. It satisfies countable chain condition if every subfamily of disjoint regions in it is countable.
\end{definition}
In a category base, a set is called a `rare set' [$7$] if its intersection with every meager set is countable. If a category base is point-meager satisfying countable chain condition and the family of all regions has cardinality at most $2^{\aleph_0}$, then under the assumption of continuum hypothesis it can be shown that every abundant set contains a rare set of cardinality $2^{\aleph_0}$ [$7$]. In a perfect base [$7$] satisfying countable chain condition, every rare set is non-Baire [$7$]. In generalizing a result of Sierpinski, Morgan showed (Th $15$, III, Ch $2$, [$7$]) that under the assumptions stated above, every abundant set is representable as the union of $2^{\aleph_1}$ non-Baire rare sets such that the intersection of any two different set  is countable. In this paper, we prove a result of a different flavour dealing with $\aleph_0$-independent family of non-Baire rare sets in any abundant Baire subset of a point-meager base satisfying countable chain condition.

\begin{definition}
	In a category base ($X,\mathcal{C}$), a set $F$ is a full subset of a set $E$ if $F\subseteq E$ and for every abundant Baire set $B$, $B\cap F$ is abundant whenever $B\cap E$ is so. If $E$ is a Baire set, this is equivalent to stating that $E-F$ cannot contain any abundant Baire set.
\end{definition}
The above definition formulated in the pattern of Grzegorek and Labuda [$2$], provides a common generalization of two analogous concepts of full subsets in measure and category within the common framework of category bases.\\

As stated above, our result in this paper involves $\aleph_0$-independent family of sets. As a particular case of the general definition given in [$10$], a family \{$A_i:i\in I$\} of sets in $X$ is called $\aleph_0$-independent (resp. strictly $\aleph_0$-independent) if for each set $J$ having card($J)<\aleph_0$ (resp. card($J)\leq \aleph_0$) and every function $f:J\mapsto \{0,1\}$, $\bigcap\{{A}_j^{f(j)}: j\in J\}\neq \emptyset $ where ${A}_j^{f(j)}=A_j$ if $f(j)=0$ and ${A}_j^{f(j)}=X-A_j$ if $f(j)=1$.\\

The existence of a $\aleph_0$-independent family of subsets of an infinite set $E$ was solved by Tarski [$6$]. He showed that such a family exists and has cardinality $2^{card(E)}$. The result has many important applications [$3$], [$10$] in solving problems related to measure extension.\\

Within the framework of category bases, the idea of a $\aleph_0$-independent (resp. strictly $\aleph_0$-independent) family can be made further strong. We say that
\begin{definition}
	A family \{$A_i:i\in I$\} of subsets of an abundant set $E$ is $\aleph_0$-independent  (resp. strictly $\aleph_0$-independent) in $E$ with respect to a category base ($X,\mathcal{C}$) if for each set $J$ with card($J)<\aleph_0$ (resp. card($J)\leq\aleph_0$ ) and every function $f\mapsto \{0,1\}$, the set $\bigcap\{{A}_j^{f(j)}: j\in J\}$ is a full subset of $E$.
\end{definition}
Let $\Phi$ be a set of one-to-one mappings of $X$ onto $X$ which is closed with respect to the composition of mappings and formation of inverses. We say that a family $\mathcal{S}$ of subsets of $X$ is invariant under $\Phi$ or simply $\Phi$-invariant [$7$] if $\phi(\mathcal{S})=\mathcal{S}$ for all $\phi\in\Phi$ where $\phi(\mathcal{S})=\{\phi(S):S\in\mathcal{S}\}$. It is easy to check that in a category base ($X,\mathcal{C}$), if $\mathcal{C}$ is a $\Phi$-invariant family, then so are the families of meager, abundant and Baire sets [$7$].
\newpage
Our theorem in this paper is as follows
\begin{theorem}
	Let ($X,\mathcal{C}$) be a point-meager, Baire base satisfying countable chain condition (ccc) and the cardinality of the class $\Phi$ is at most $2^{\aleph_0}$. Moreover, let $\mathcal{C}$ be $\Phi$-invariant having cardinality at most $2^{\aleph_0}$ and satisfying the condition that for any two regions $C$ and $D$, there exists $\phi\in\Phi$ such that $C\cap\phi(D)$ contains a region. Then under continuum hypothesis in every abundant Baire set $E$ satisfying $\phi(E)= E$ for all $\phi\in\Phi$, there exists a family of non-Baire rare sets having cardinality $2^{\aleph_1}$ such that for every set $A$ in this family and every $\phi\in\Phi$, $\phi(A)$ is also a rare set and this family is strictly $\aleph_0$-independent and hence $\aleph_0$-independent in $E$ with respect to ($X,\mathcal{C}$).
\end{theorem}
A proof of the above theorem depends on the following lemma which is a particular case of Proposition $2.10$ stated in [$10$].
\begin{lemma}
If $E$ is an infinite set satisfying  the condition (card $E)^{\aleph_0}=$card $E$, then there exists a maximal strictly $\aleph_0$-independent family \{$A_i:i\in I$\} of subsets of $E$ such that card($I)=2^{card(E)}$.
\end{lemma}	 
In our proof, we also utilize the concept of an Ulam ($\aleph_0$,$\aleph_1$)-matrix [$5$] the definition of which is given below.
\begin{definition}
Let $E$ be an infinite set with card($E)=\aleph_1$. A double family ($E_{\xi,\zeta})
_{\xi<\aleph_0,\zeta<\aleph_1}$ of subsets of $E$ is called an Ulam ($\aleph_0$,$\aleph_1$)-matrix over $E$ if the following two conditions are satisfied:
\begin{enumerate}
	\item card($E-\bigcup\{E_{\xi,\zeta}:\xi<\aleph_0\})\leq\aleph_0$ for every $\zeta<\aleph_1$.
	\item $E_{\xi,\zeta}\cap E_{\xi,\zeta^\prime}=\emptyset$ for all $\xi<\aleph_0$ and any two distinct ordinals $\zeta<\aleph_1,\zeta^\prime<\aleph_1$.
\end{enumerate}
\end{definition}

 \begin{proof} [Proof of  Theorem 2.6]
Since $\mathcal{C}$ satisfies countable chain condition having cardinality at most $2^{\aleph_0}$, every meager set in ($X,\mathcal{C})$ is contained in a $\mathcal{K}_{\delta\sigma}$-meager set (Th $5$, II, Ch $1$, [$7$]) where $A\in\mathcal{K}$ iff $X-A\in\mathcal{C}$ and moreover as $X$ is abundant, card($\mathcal{K}_{\delta\sigma})=2^{\aleph_0}$ (Th $1$, I, Ch $2$, [$7$]). Let $\Omega$ be the smallest ordinal representing card($E$) (we shall use the same notations $\aleph_0$, $\aleph_1$ for the smallest ordinals representing the cardinals $\aleph_0$ and $\aleph_1$) and consider the well orderings
\begin{equation*}
\begin{matrix}
	x_1, & x_2, & \cdots & x_\alpha, & \cdots (\alpha<\Omega)\\
	R_1, & R_2, & \cdots & R_\alpha, & \cdots (\alpha<\aleph_1)\\
	\phi_1, & \phi_2, & \cdots & \phi_\alpha, & \cdots (\alpha<\aleph_1)
\end{matrix}
\end{equation*}
of $E$, the family of all $\mathcal{K}_{\delta\sigma}$-meager sets and the family $\Phi$ of all mappings. These arrangements are justified, because according to the continuum hypothesis, we have $2^{\aleph_0}=\aleph_1$.\\
We now proceed to construct a subset of $E$ in the following manner by defining a transfinite sequence 
\begin{equation*}
	y_1, ~ y_2, ~ \cdots ~ y_\alpha,  ~\cdots~ (\alpha<\aleph_1)
\end{equation*} 
Select arbitrarily a point $x_1$ from $E$ and put $y_1=x_1$. Assume that for any ordinal $\alpha$ ($\alpha<\aleph_1$), we have already selected elements $y_1$, $y_2$, $\cdots$, $y_\beta$ $(\beta<\alpha)$. Let $G_\alpha$ be the group generated by the mappings $\phi_\beta$ ($\beta<\alpha$), i.e. $G_\alpha$ consists of all the elements which are of the form ${\phi}_{\alpha_1}^{n_1}*{\phi}_{\alpha_2}^{n_2}*\cdots*{\phi}_{\alpha_k}^{n_k}$ where $k\in\mathbb{N}$ (the set of naturals); $n_1$, $n_2$, $\cdots$, $n_k$ are integers and ordinals $\alpha_1$, $\alpha_2$, $\cdots$, $\alpha_k$ ($<\alpha$).\\
We set $P_\alpha=\bigcup\limits_{\beta<\alpha}\{\psi^{-1}(R_\beta):\psi\in G_\alpha\}$ and $B_\alpha=\{\psi(y_\beta):\beta<\alpha, \psi\in G_\alpha\}$. The set $B_\alpha$ is countable and hence meager by hypothesis. Also $P_\alpha$ is meager because $\mathcal{C}$ is a $\Phi$-invariant family. Hence $P_\alpha\cup B_\alpha$ is meager. We choose  $y_\alpha$ as the first element from $E$ which lies outside $P_\alpha\cup B_\alpha$.\\
The sets $B_\alpha$ ($1<\alpha<\aleph_1$) along with $B_1=\emptyset$ forms an increasing transfinite sequence, i.e. $B_\beta\subseteq B_\alpha$ if $\beta<\alpha<\aleph_1$ and $B_\alpha=\bigcup\limits_{\beta<\alpha}{B_\beta}$ if $\alpha$ is a limit ordinal. We set $O_\alpha=B_{\alpha+1}-B_\alpha$. Since $y_\alpha=\psi^{(0)}y_\alpha$, so $y_\alpha\in O_\alpha$ and hence $O_\alpha\neq\emptyset$ for every $\alpha$. Moreover, the sets $O_\alpha$ are mutual disjoint, i.e. $O_\alpha\cap O_\beta=\emptyset$ if $\alpha\neq\beta$ by construction. We set $A=\bigcup\limits_{\alpha<\aleph_1}{O_\alpha}$ which is obviously an abundant set contained in $E$.\\
Now we consider the Ulam ($\aleph_0$,$\aleph_1$)-matrix ${({\Pi}_{\xi,\zeta})}_{\xi<\aleph_0,\zeta<\aleph_1}$ on $\aleph_1$ and set ${E}_{\xi,\zeta}=\bigcup\limits_{\gamma\in{\Pi}_{\xi,\zeta}}O_\gamma$. Then there exists $\xi_0$ and a subset $\Theta$ of $\aleph_1$ having card($\Theta)=\aleph_1$ such that ${E}_{\xi_0,\zeta}$ is abundant in ($X,\mathcal{C}$) for every $\zeta\in\Theta$ and for any two $\zeta$, $\zeta^\prime$ ($\zeta\neq\zeta^\prime$), ${E}_{\xi_0,\zeta}\cap{E}_{\xi_0,\zeta^\prime}=\emptyset$. This is so because ($X,\mathcal{C}$) is point-meager, $A$ is abundant and the class of meager sets in ($X,\mathcal{C}$) forms a $\sigma$-ideal. By continuum hypothesis, $2^{\aleph_0}=\aleph_1$ and therefore ${\aleph_1}^{\aleph_0}={(2^{\aleph_0})}^{\aleph_0}=2^{\aleph_0}=\aleph_1$. So by Lemma $2.7$, there exists a strictly $\aleph_0$-independent family \{$\Theta_i:i\in I$\} of subsets of 
$\Theta$ such that card($I)=2^{\aleph_1}$. Consequently, for every set $J\subseteq I$ having card($J)\leq\aleph_0$ and every function $f:J\mapsto\{0,1\}$,  $\bigcap\limits_{j\in J}{A}_j^{f(j)}\neq\emptyset$ where $A_i=\bigcup\limits_{\zeta\in\Theta_i}E_{\xi_0,\zeta}$.\\
Evidently, from the construction it follows that each $A_i$ is a rare set. Also for every $\psi\in\Phi$ and every $A_i$ in the family \{$A_i: i\in I$\}, $\psi(A_i)\triangle A_i$ is countable and therefore $\psi(A_i)$ is a rare set. We claim that no $A_i$ can be a Baire set by showing that every abundant Baire set meets any set $T$ which is the union of some $E_{\xi_{0},\zeta}$ non-vacously. If possible, let there be an abundant Baire set $B$ such that $B\cap T=\emptyset$. Let $B$ be abundant everywhere in a region $C$ and $T$ be abundant everywhere in a region $D$. By hypothesis, there exists $\phi^*\in\Phi$ such that $C\cap\phi^*(D)$ contains a region $F$ (say). Consequently, $F-B$ is meager (Th 2, III, Ch 1, [7]) and ($X,\mathcal{C}$) being a Baire base, $B\cap\phi^*(T)$ is abundant. But this is impossible because $B$ is disjoint from $T$ and $\phi^*(T)\triangle T$ is countable and hence meager. Therefore, $B\cap T\neq \emptyset$. But in the above construction we can always choose $B=A_i$ and $T=E_{\xi_{0},\zeta_i}$ such that $A_i\cap E_{\xi_0,\zeta_i}=\emptyset$ even when $A_i$ is a Baire set.\\
Lastly, the family \{$A_i:i\in I$\} is strictly $\aleph_0$-independent with respect to ($X,\mathcal{C}$). We set $T=\bigcap\limits_{j\in J}{A}_j^{f(j)}$ where $J\subseteq I$ with card($J)\leq\aleph_0$ and $f:J\mapsto\{0,1\}$ as above. Since every set $T$ can be expressed as $T=\bigcup\{E_{\xi_0,\zeta}:\zeta\in \bigcap\limits_{j\in J}{\Theta}_j^{f(j)} \}$ where $\bigcap\limits_{j\in J}{\Theta}_j^{f(j)}\neq\emptyset$ and each $E_{\xi_0,\zeta}$ is abundant, so by the same reasoning as given above, $T$ meets every abundant Baire set contained in $E$.
\end{proof}
\begin{remark}
	Let $X=\mathbb{R}$ and $\Phi$ represents the group of all translations in $\mathbb{R}$. Then for the category base ($\mathbb{R},\mathcal{C}$) where $\mathcal{C}$ is the usual topology of $\mathbb{R}$, the hypothesis of Theorem 2.6 is satisfied and in $\mathbb{R}$ we can have a family of Luzin sets without Baire property such that for every set $A$ in this family and every $x\in\mathbb{R}$, $A+x$ is also a Luzin set which obviously follows from the translation invariance of first category sets and the family is strictly $\aleph_0$-independent with respect to the usual topology of $\mathbb{R}$. As it is also true (by regularity of Lebesgue measure) that for every measurable set $A$ of positive measure, $A\cap(A+x)$ is a set of positive measure for every $x$ in some interval around the origin, so for the category base ($\mathbb{R},\mathcal{C}$) where $\mathcal{C}$ is the family of Lebesgue measurable sets of positive measures, the hypothesis of Theorem 2.6 is satisfied, and so in $\mathbb{R}$ we can have a family of nonmeasurable Sierpinski sets such that for every set $A$ in this family and every $x\in\mathbb{R}$, $A+x$ is also a Sierpinski set  which obviously follows from the translation invariance of measure zero sets and the family is strictly $\aleph_0$-independent with respect to the Lebesgue measure space.
	\end{remark}

\bibliographystyle{plain}

	\end{document}